\documentclass[12pt]{amsart}
\usepackage{amssymb}
\usepackage{amsthm}
\usepackage{amsmath}
\usepackage{color}
\usepackage{float}
\usepackage{graphicx}
\usepackage[top=1in, bottom=1in, left=1in, right=1in]{geometry}

\newtheorem{Definition}{Definition}
\newtheorem{Example}{Example}[section]
\newtheorem{Theorem}{Theorem}[section]
\newtheorem{Lemma}{Lemma}[section]
\newtheorem{Corollary}{Corollary}[section]
\newtheorem{Conjecture}{Conjecture}
\newtheorem*{remark}{Remark}
\newtheorem*{Remark}{Remark}
\newcommand{\bigslant}[2]{{\raisebox{.2em}{$#1$}\left/\raisebox{-.2em}{$#2$}\right.}}

\newcommand{\I}{\mathfrak{B}I}

\title{Diagrammatic Calculus of Coxeter and braid groups} 
\author{Niket Gowravaram and Uma Roy \\ Mentor: Alisa Knizel \\MIT PRIMES}

\begin{document}

\begin{abstract}
We investigate a novel diagrammatic approach to examining strict actions of a Coxeter group or a braid group on a category. This diagrammatic language, which was developed in a series of papers by Elias, Khovanov and Williamson, provides new tools and methods to attack many problems of current interest in representation theory. In our research we considered a particular problem which arises in this context. To a Coxeter group $W$ one can associate a real hyperplane arrangement, and can consider the complement of these hyperplanes in the complexification $Y_W$. The celebrated $K(\pi,1)$ conjecture states that $Y_W$ should be a classifying space for the pure braid group, and thus a natural quotient $\bigslant{Y_W}{W}$ should be a classifying space for the braid group. Salvetti provided a cell complex realization of the quotient, which we refer to as the Salvetti complex. In this paper we investigate a part of the $K(\pi,1)$ conjecture, which we call the $K(\pi,1)$ conjecturette, that states that the second homotopy group of the Salvetti complex is trivial. In this paper we present a diagrammatic proof of the $K(\pi,1)$ conjecturette for a family of braid groups as well as an analogous result for several families of Coxeter groups. 
\end{abstract}

\maketitle

\section{Introduction} 

Group theory, which is the study of algebraic structures known as groups, is a vitally important part of mathematics that has applications in various fields including physics, chemistry, crystallography, cryptography, and combinatorics, as well as being a rich area of study in its own right \cite{crypto, Noether, crystals, quantum, chem, Hum}. Two groups that arise often in the study of natural phenomenon are the dihedral group $D_n$ and the symmetric group $S_n$. The dihedral group and symmetric group are both special cases of a more general class of groups known as Coxeter groups---the main focus of our project.

In addition to being generalizations of the natural reflection groups, Coxeter groups have a myriad of uses in mathematics, especially in representation theory, where they serve as building blocks for the classification of algebraic objects. Examples of finite Coxeter groups include the symmetry groups of regular polytopes and the Weyl groups of simple Lie algebras, which are very important in the study of particle physics \cite{quantum}. Infinite Coxeter groups include symmetry groups of regular tessellations of Euclidean space and Weyl groups of affine Kac-Moody algebras, which are a generalization of semisimple Lie algebras and are of particular importance in conformal field theory and the theory of exactly solvable models \cite{Hum, quantum}.
    
To a Coxeter group $W$ one can associate a real hyperplane arrangement and consider the complement of these hyperplanes in the complexification $Y_W$. The celebrated $K(\pi,1)$ conjecture in modern algebraic topology states that $Y_W$ should be a classifying space for the pure braid group, and thus a natural quotient $\bigslant{Y_W}{W}$ should be a classifying space for the braid group. In \cite{salvetti}, Salvetti provided a cell complex realization of the quotient, which we refer to as the \emph{Salvetti complex}. The $K(\pi,1)$ conjecture was proven for finite Coxeter groups by Deligne in \cite{deligne} but many cases remain open. In this paper we use a novel approach to investigate a part of the $K(\pi,1)$ conjecture, which we refer to as the $K(\pi,1)$ conjecturette, that states that the second homotopy group (denoted as $\pi_2$) of the Salvetti complex is trivial. In \cite{elias}, another cell complex was introduced as a $3$-skeletal model for the classifying space of a Coxeter group $W$. In this paper, we also prove that $\pi_2$ of this cell complex is trivial for several series of finite Coxeter groups, verifying that it is indeed a valid $3$-skeletal approximation. 

Due to a diagrammatic interpretation of $\pi_2$ which can be found in a book by Fenn \cite{fenn} one can think of the elements of $\pi_2$ of the Salvetti and the cell complexes introduced in \cite{elias} as special types of decorated planar graphs, which we refer to as \emph{diagrams}. Two diagrams are considered homotopic if one can be transformed into the other using a sequence of allowed transformations, which we  describe in Section \ref{background} of our paper. The problem we are considering naturally splits into two directions. One is studying unoriented diagrams, which corresponds to Coxeter groups and the topology of their classifying spaces, and the other is studying similar diagrams but with orientations on the edges, which corresponds to braid groups and the topology of the Salvetti complex.  The goal of our project was to prove that any diagram is homotopic to the empty diagram in both cases (Conjecture~\ref{main}), which is equivalent to the triviality of $\pi_2$ of the corresponding complexes (more details can be found in \cite{elias}). The beauty of this diagrammatic approach is the elementary nature of our combinatorial methods, which are used to prove deep statements in modern algebraic topology.
    
In this paper, we present our results for the symmetric groups and dihedral groups. We also present our results for the Artin braid group $\I_n$, which is a generalization of the dihedral Coxeter group. Our diagrammatic proof of the $K(\pi,1)$ conjecturette for the aforementioned family of braid groups and our results on Coxeter groups answer a question posed in \cite{elias} regarding the existence of diagrammatic proofs for these type of statements. Our findings represent research towards proving the $K(\pi,1)$ conjecturette for all Coxeter groups diagrammatically. In addition, our proof for the braid group $\I_m$ is among the first proofs, to our knowledge, using the diagrammatic calculus for braid groups developed in \cite{elias}. Our paper is organized as follows: In the background section we introduce all of the definitions needed for our work. In Section \ref{prelims} we state several general theorems and lemmas, and as a consequence derive our result for the dihedral groups. In Sections \ref{secAn} and \ref{In} we prove Conjecture ~\ref{main} for the aforementioned families of Coxeter and braid groups. We conclude our work in Section \ref{conc}.

\section{Background} \label{background} 

\subsection{Coxeter Groups}
We begin by introducing some basic definitions associated with the study of Coxeter groups.

\begin{Definition}[Coxeter Group] 

A Coxeter group is a group given by generators $g_1, \ldots, g_n$ with relations $(g_ig_j)^{m_{i,j}}=1$ for each pair of generators $g_i, g_j$,
such that $m_{i,j} \in \mathbb{N}$ where $m_{i,j} \geq 2$ for $i \neq j$ and $m_{i,i} = 1$ for all $i$.

\end{Definition}

\begin{remark}
Some numbers $m_{i,j}$ can be $\infty$, in which case there is no relation between generators $g_i$ and $g_j$. The condition $m_{i,i} = 1$ implies  $g_i^2=1$ and as a consequence $m_{i,j}=m_{j,i}$.

\end{remark}

\begin{Example}\label{symgroup}
The symmetric group $S_n$ has generators $g_1, g_2, \ldots, g_{n-1}$, where $g_i$ is $i^{th}$ elementary transposition that sends $i \rightarrow i+1$ and $i+1 \rightarrow i$, subject to the following relations$\colon$ $g_i^2=1$, $(g_ig_{i+1})^3 = 1$, and $(g_ig_j)^2 = 1$ if $j \neq i+1, i-1$.
\end{Example}

To each Coxeter group, we can associate a Coxeter-Dynkin diagram, which are a convenient method of visualizing the generators and relations of a Coxeter group and also useful in classifying the Coxeter groups.

\begin{Definition}
For a particular Coxeter group, the associated Coxeter-Dynkin diagram is a graph where vertices correspond to generators of the group, and edges correspond to relations between generators. We let $v_i$ and $v_j$ be arbitrary vertices in our diagram corresponding to the generators $g_i$ and $g_j$ respectively. We have the following 3 properties:

\begin{itemize}
\item If there is no edge between $v_i$ and $v_j$, then $g_i$ and $g_j$ commute (i.e. $m_{i,j} = 2)$.
\item If there is an unlabelled edge between $v_i$ and $v_j$, then $m_{i,j} = 3$. 
\item If there is an edge labeled with an integer $k$, for some $k \in \mathbb{N}$, between $v_i$ and $v_j$ then $m_{i,j} = k$. 
\end{itemize}

\end{Definition}

\begin{Example}
The symmetric group $S_n$ has the following Coxeter-Dynkin diagram with $n-1$ vertices:

\end{Example}

\begin{figure}[H]
\centering
\includegraphics[scale = 0.3]{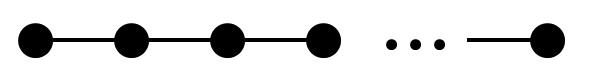} 
\caption{The Coxeter-Dynkin diagram of the symmetric group.}
\end{figure}

\emph{Irreducible Coxeter groups} are Coxeter groups that have connected Coxeter-Dynkin diagrams. All Coxeter groups are the direct product of irreducible Coxeter groups. In our research, we only need to consider the irreducible Coxeter groups due to Theorem ~\ref{GcrossH}. Using Coxeter-Dynkin diagrams, one can classify all finite, irreducible Coxeter groups \cite{Hum}. We note that the Coxeter group $A_n$ is isomorphic to the symmetric group $S_n$ and the Coxeter group $I_n$ is isomorphic to the dihedral group $D_n$. For the rest of this paper, we use $A_n$ and $I_n$ in place of $S_n$ and $D_n$ respectively.

\subsection{Diagrammatics of Coxeter groups}
Given a Coxeter group $W$ we assign each generator a unique color and color the vertices of the Coxeter-Dynkin diagram by the color of the corresponding generator. Having assigned each generator a color, we consider certain diagrams, which are colored planar graphs along with a number of distinct circles (edges with no attached vertices). Every edge of the diagram is colored by a color corresponding to one of the generators of the Coxeter group. Every vertex of the diagram must correspond to a pair of generators $g_i$ and $g_j$, and must have degree $2m_{i,j}$ with edges alternating between the colors of the two generators. Note that there can be many vertices corresponding to the same pair of generators.

\indent Two diagrams are considered homotopic if one can be transformed into another through a series of the following transformations: 

    
\noindent  \textbf{Circle relation: } Given a graph, we can add or remove empty circles of any color.

    
    
\noindent  \textbf{Bridge relation: } Given two edges of the same color, we can switch around the vertices they connect to as long as we do not create any new intersections.
    
\begin{figure}[H]
\centering
\includegraphics[scale = 0.2]{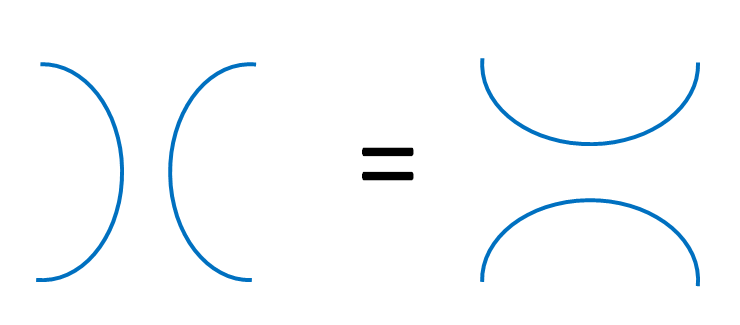} 
\caption{A drawing of the bridge relation.}
\end{figure}


\begin{remark}
Note that the circle and bridge relation allow us to remove all circles in our graph. Thus we assume we have no circles in our diagrams from this point forth.
\end{remark}

\noindent \textbf{Cancellation of pairs relation: } For any subgroup of a Coxeter group $G$ isomorphic to $I_m$ we have the following allowed transformation. Note that that we draw the relation for $I_3$, but the relation holds for general $I_m$.

\begin{figure}[H]
\centering
\includegraphics[scale = 0.2]{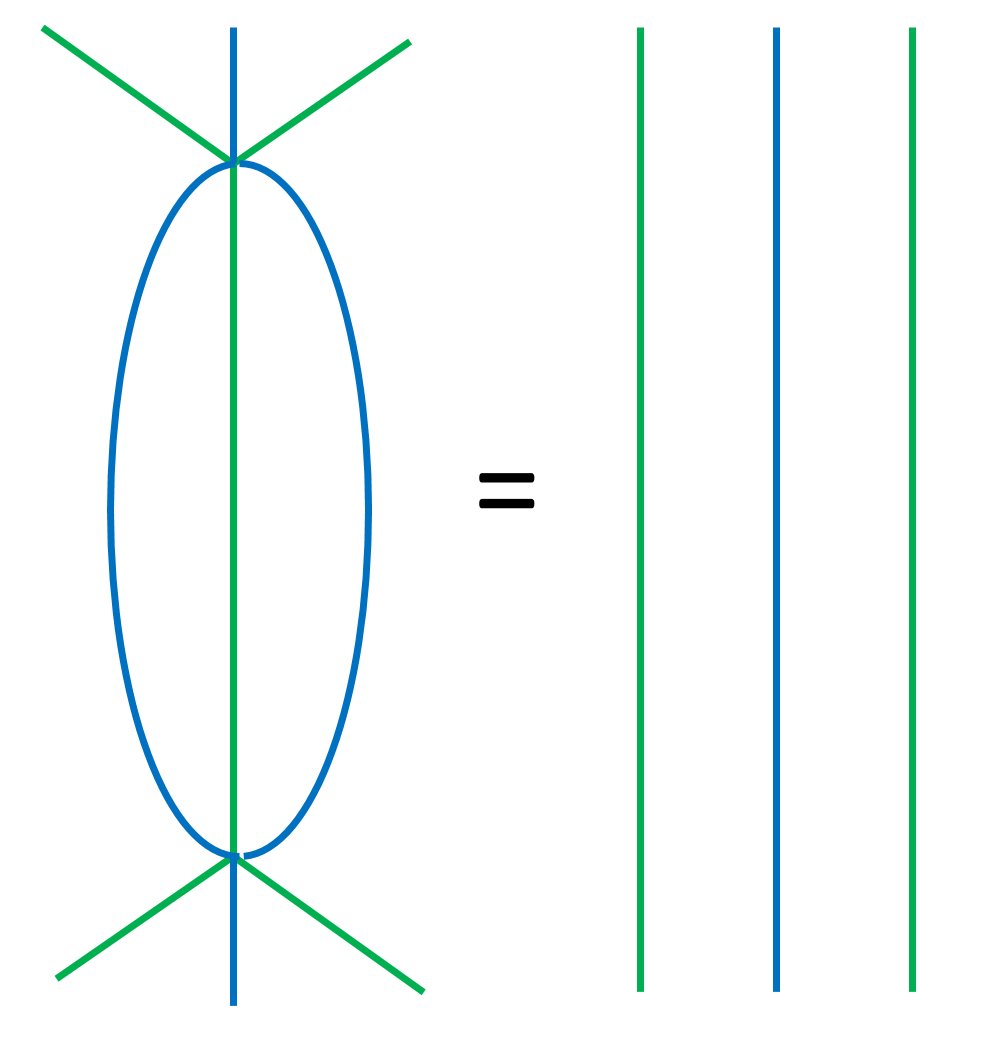} 
\caption{The cancellation of pairs relation corresponding to $I_3$. In general, the cancellation of pairs relation for $I_m$ looks similar to the above picture, however each vertex has degree $2m$.}
\label{Imzam}
\end{figure}

Finally, we have a class of relations called the \emph{Zamolodchikov} relations (written as ZAM relations for brevity) that are determined by inspection of the reduced expression graph for the longest element of the finite rank 3 Coxeter groups: $A_3, B_3, H_3$ and $A_1\times I_n$ (more detail can be found in \cite{elias}). All of the ZAM relations we use in our paper are drawn below. Note that while we draw the ZAM relations below with specific colors, they hold true for any colors that form the arrangement of vertices and edges in the drawings below.

\begin{figure}[H]
\centering
\includegraphics[scale = 0.3]{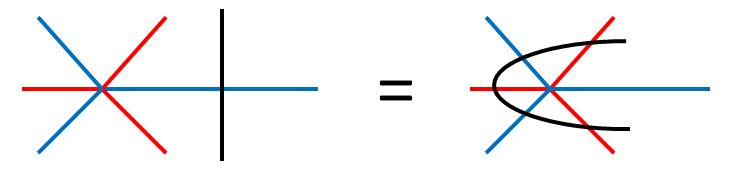} 
\caption{The ZAM relation corresponding to the group $A_1\times A_2$.}
\label{flipzam}
\end{figure}

\begin{figure}[H]
\centering
\includegraphics[scale = 0.2]{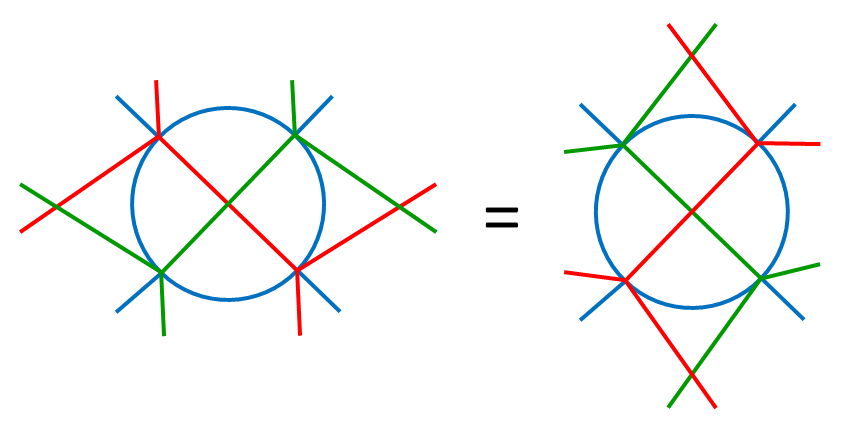} 
\caption{The ZAM relation corresponding to the group $A_3$.}
\label{a3zam}
\end{figure}

\subsection{Braid groups and their diagrammatics}

We also consider the diagrammatics of Artin braid groups, which are a generalization of Coxeter groups. The Artin braid groups are far more complicated to work with diagrammatically, as our diagrams become oriented planar graphs with a number of oriented circles (as we will describe below).

\begin{Definition}
An Artin braid group is given by generators $g_1, \ldots, g_n$ with relations $(g_ig_j)^{m_{i,j}}=1$ between every pair of generators $g_i, g_j$ with $i \neq j$ such that $m_{i,j} \in \mathbb{N}$ and $m_{i,j} \geq 2$. There is no longer any relation $m_{i,i}$ as there was in the Coxeter groups.

\end{Definition}

In this paper we only work with the braid group $\I_m$, which is defined below.

\begin{Definition}\label{braidIn}

The braid group $\I_m$ is given by 2 generators $g_1$ and $g_2$ with $m_{1,2} = m$. For the rest of our paper, when considering $g_1$ and $g_2$ as part of $\I_m$, we assign $g_1$ the color blue and $g_2$ the color green.

\end{Definition}

Similar to how we constructed diagrams for Coxeter groups, we do the same for the braid groups as follows. For a given braid group $\mathfrak{B}$ with generators $g_1, g_2, \ldots, g_n$ and with relations $(g_ig_j)^{m_{i,j}}=1$ between generators, we assign a distinct color to each generator. Every edge of the graph must be a color corresponding to a generator. Every vertex of the graph must correspond to a pair of generators, $g_i$ and $g_j$, and must have degree $2m_{i,j}$ with edges alternating between the colors of the two generators. In addition, our edges have orientations as specified: each vertex must have $m_{i,j}$ consecutive edges pointing out of the vertex and $m_{i,j}$ consecutive edges pointing towards the vertex. One can see that orienting edges in this manner results in two distinct types of vertices for each pair of generators.
\vspace{-2mm}

\begin{Example}
For the braid group $\I_3$, there are 2 different types of vertices, drawn below in Figure~\ref{2diffverts}. Recall that we have assigned $g_1$ the color blue and $g_2$ the color green.

\end{Example}

\begin{figure}[H]
\centering
\includegraphics[scale = 0.5]{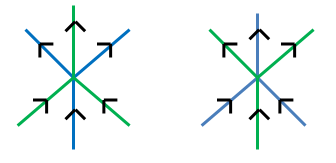}
\caption{The 2 different types of vertices in diagrams for the braid group $\I_3$. }
\label{2diffverts}
\end{figure}

Again, two diagrams are considered homotopic if one can be transformed into the other through the following 3 transformations. Note that these are the same transformations as for the diagrammatics of Coxeter groups with added orientations.

    
\noindent \textbf{Circle relation: } Given a graph, we can add or remove empty oriented circles of any color.

We note that as for the Coxeter circle relation, this relation allows us to ignore any oriented circles in our diagrams, as we can simply remove them.
    
    

\noindent  \textbf{Directed bridge relation: } Given two edges of the same color, we can switch around the vertices they connect to as long as we do not create any new intersections and the orientations of the edges are preserved.

\begin{figure}[H]
\centering
\includegraphics[scale = 0.2]{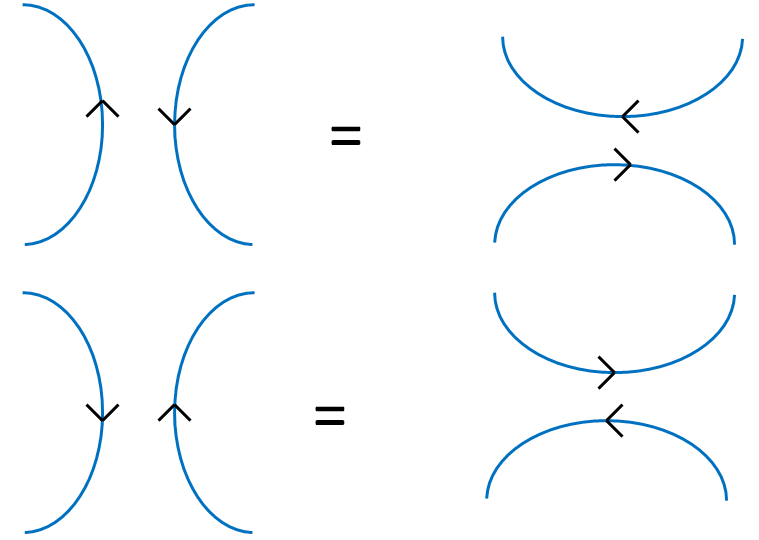} 
\caption{A drawing of the directed bridge relation.}
\end{figure}



\noindent  \textbf{Directed cancellation of pairs: } For any subgroup of a braid group $\mathfrak{B}$ isomorphic to $\I_m$ we have the following allowed transformation. Note that that we draw the relation for the group $\I_3$ but it holds for general $\I_m$.

\begin{figure}[H]
\centering
\includegraphics[scale = 0.4]{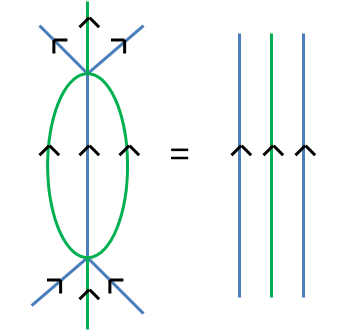} 
\caption{The cancellation of pairs relation corresponding to the group $\I_3$.}
\label{directedImzam}
\end{figure}

In general, the cancellation of pairs relation for the group $\I_m$ looks similar to Figure ~\ref{directedImzam}, however each vertex has degree $2m$. This cancellation of pairs relation also holds if all arrows in the above picture are reversed in orientation (i.e. they all point downwards).
The braid groups also have ZAM relations corresponding to rank $3$ subgroups, however we will not need any in our proofs and thus do not include them in this paper. 

\subsection{The $K(\pi,1)$ conjecturette}

The goal of our project was to prove the following conjecture.

\begin{Conjecture}[$K(\pi,1)$ conjecturette]\label{main}
All diagrams for a particular Coxeter group or braid group are homotopic to the empty diagram through a series of the aforementioned allowed transformations.

\end{Conjecture}

\begin{remark}

As mentioned in our introduction, Conjecture ~\ref{main} for braid groups is equivalent to showing that $\pi_2$ of the Salvetti complex of the corresponding Coxeter group is trivial, which is a part of the $K(\pi,1)$ conjecture in representation theory. Conjecture ~\ref{main} for Coxeter groups has a similar topological interpretation (see \cite{elias}), which is equivalent to showing that $\pi_2$ of the cell complex introduced in \cite{elias} is trivial, verifying that the cell complex is a valid 3-skeletal approximation for the classifying space of a Coxeter group $W$. We remark that Conjecture ~\ref{main} for Coxeter groups is not part of the $K(\pi,1)$ conjecture, as the $K(\pi,1)$ conjecture deals with the topology of Salvetti complexes and Conjecture ~\ref{main} for Coxeter groups deals with the cell complex introduced in \cite{elias}. However, for the purposes of this paper, we refer to Conjecture ~\ref{main} for \emph{both} Coxeter and braid groups as the $K(\pi,1)$ conjecturette. In \cite{elias}, authors Elias and Williamson raised the question of the existence of an elementary diagrammatic proof for Conjecture ~\ref{main} for both Coxeter and braid groups. The results in our paper address this question by providing diagrammatic proofs for Conjecture ~\ref{main} for several families of Coxeter groups and a family of braid groups. 

\end{remark}

\section{General Statements}\label{prelims}

In this section we provide several general statements that are useful in many of our proofs. 

\begin{Theorem}\label{GcrossH}
Given Coxeter groups $G$ and $H$, if Conjecture ~\ref{main} holds for $G$ and for $H$, then it holds for the group $G \times H$.

\end{Theorem}

\begin{remark}
Due to the above theorem, we see that if we prove Conjecture ~\ref{main} for all irreducible Coxeter groups, then we have proven it for all Coxeter groups, as all Coxeter groups are the direct product of the irreducible Coxeter groups. This allows to consider only the irreducible Coxeter groups in our research. We note that the above theorem also holds if $G$ and $H$ are braid groups. The proof when $G$ and $H$ are braid groups is the same as the proof for when $G$ and $H$ are Coxeter groups, except it uses the braid version of the ZAM relation for $A_1 \times I_m$.
\end{remark}

\begin{proof}
We begin by observing that any generator $g \in G$ commutes with any generator $h \in H$. Given a diagram, consider the subgraph with only edges corresponding to generators of $G$. We call this subgraph the \emph{$g$-subgraph} of the diagram. Examining any 2 adjacent vertices on the $g$-subgraph connected by an edge $E$, we see that edge $E$ may be intersected by edges corresponding to generators $h \in H$. Thus using the ZAM relation corresponding to $A_1 \times I_m$ (Figure ~\ref{flipzam}), since any pair of generators $g,h$ with $g \in G$ and $h \in H$ commute, we can remove all edges intersecting $E$. By continuing this process on all edges of the $g$-subgraph, we can remove all edges corresponding to generators in $H$ from the $g$-subgraph. Thus we end up with 2 disjoint graphs---one that has only edges corresponding to generators in $G$ and the other that has only edges corresponding to generators in $H$. Since we know all diagrams for groups $G$ and $H$ are homotopic to the empty graph, we can reduce these 2 graphs to the empty graph. Thus, all possible diagrams of $G \times H$ are homotopic to the empty one. 
\end{proof}


\begin{Lemma}\label{adjvert}
If two vertices of the same type (vertices that have edges of the same two colors) are connected by an edge, then we can delete both the vertices.
\end{Lemma}

\begin{proof}
Consider two connected vertices of the same type where each vertex has degree $2m$. Since these two vertices are connected by an edge, we can use the bridge relation locally to connect the other edges of these two vertices. Using the cancellation of pairs relation for $I_m$ (Figure ~\ref{Imzam}), we can remove both of these vertices. 
\end{proof}

Below is a picture showing this process for the Coxeter group $I_3$ (which is equivalently $A_2$, the symmetric group of order 2). 

\begin{figure}[H]
\centering
\includegraphics[scale = 0.16]{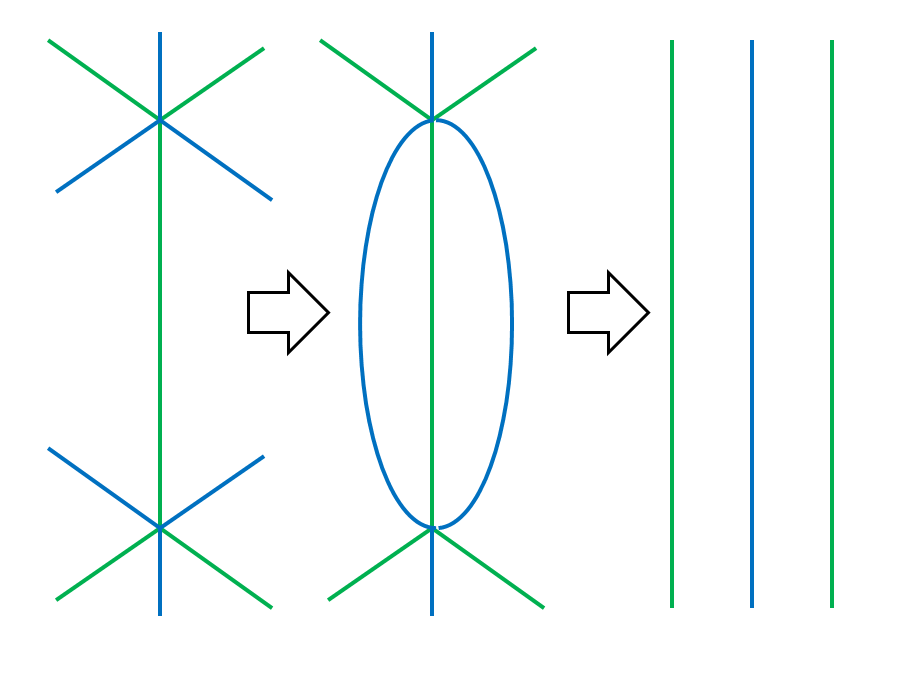} 

\caption{Deleting 2 adjacent vertices of the same type in the Coxeter group $I_3$. The first step comes from using the bridge relation. The second step comes from using the $I_3$ cancellation of pairs relation.}
\end{figure}

\begin{Corollary}\label{cordihedral}
Conjecture ~\ref{main} is true for all Coxeter groups $I_n$. The Coxeter groups $I_n$ are isomorphic to the dihedral groups. 
\end{Corollary}

\begin{proof}
Given the Coxeter group $I_n$, we see from its definition that there are only 2 generators. Thus, there is only 1 type of vertex we can form, and thus all possible diagrams will always have connected vertices of the same type. Also, as can be seen by inspection, we cannot have a vertex connect to itself or just one vertex in a connected component of the graph. Using Lemma ~\ref{adjvert}, we can delete these adjacent vertices, until there are no vertices left in our diagram and it is the empty diagram. A
\end{proof}

Before introducing the next 2 lemmas, we define \emph{boundary} and \emph{subdiagram}---two terms that help us study local properties of diagrams.

\begin{Definition}
A \emph{subdiagram} of a diagram is a subset of vertices and edges of the diagram that are connected.
\end{Definition}

\begin{Definition}
Given a subdiagram, we call edges that have exactly 1 endpoint in our subdiagram and the other endpoint a vertex not in our subdiagram the \emph{boundary} of our subdiagram.  
\end{Definition}

\begin{figure}[H]
\centering
\includegraphics[scale = 0.2]{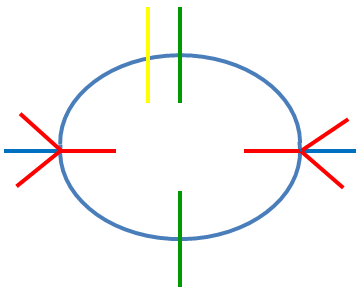}
\caption{This subdiagram has boundary yellow-green-red-blue-red-green-red-blue-red.}
\label{subdiagfig}
\end{figure}

We note that the entire diagram has an empty boundary.
\begin{Example}
Figure~\ref{subdiagfig} provides an example of a boundary.
\end{Example}

\begin{Definition}
Given 2 subdiagrams $D_1$ and $D_2$ with the same boundary, we let $D_1 \sqcup D_2$ denote the diagram obtained by connecting all edges on the boundary of $D_1$ to the corresponding edges on the boundary of $D_2$. 
\end{Definition}

\begin{Definition}
2 subdiagrams are referred to as \emph{equivalent} if they are homotopic. 

\end{Definition}

\begin{Lemma}\label{boundarycor}
If $D_1$ and $D_2$ are two subdiagrams with the same boundary and the diagram $D_1 \sqcup D_2$ is homotopic to the empty diagram (denoted as $\varnothing$), then $D_1$ is homotopic to $D_2$.
\end{Lemma}

\begin{proof}
Clearly $D_1$ is equal to the diagram with $D_1$ and the empty graph next to it. Since $D_1 \sqcup D_2 = \varnothing$, we see that $D_1$ is equivalent to $D_1$ with the closed diagram $D_1 \sqcup D_2$ next to it. Using the bridge relation on this diagram to connect all the edges belonging to $D_1$ in $D_1 \sqcup D_2$ to $D_1$, we see that we are left with $D_2$ and $D_1 \sqcup D_1$. But $D_1 \sqcup D_1 = \varnothing$ since all vertices that are connected to each other are of the same type, and thus we can delete them. Thus we are left with $D_2$. Therefore we conclude that $D_1$ is homotopic to $D_2$. 
\end{proof}

\begin{Corollary}\label{boundary}
If Conjecture ~\ref{main} holds for a Coxeter group $W$ then any 2 subdiagrams $D_1$ and $D_2$ for $W$ with the same boundary are equivalent.
\end{Corollary}

\begin{proof}

Since $D_1$ and $D_2$ have the same boundary, $D_1 \sqcup D_2$ is a closed diagram and thus is homotopic to the empty graph, since we know all diagrams for $W$ satisfy Conjecture ~\ref{main}. Thus by Lemma ~\ref{boundarycor}, we see that these 2 subdiagrams are equivalent. 
\end{proof}

\begin{Lemma}\label{trivialboundary}
When the boundary of a subdiagram is written as a word of the Coxeter group, the word is equivalent to the trivial word.
\end{Lemma}

\begin{proof} 
Each edge in a diagram represents a specific generator (as we color the edges for this purpose). Thus, the boundary of a subdiagram represents a word formed by the generators of the group. Also, we see that every vertex in our diagram represents an equivalent rewriting of a word, because of our group relations. Thus the word formed by one part of the boundary region is transformed, and thereby, equivalent to the word formed by the other part of the boundary region through all vertices in the diagram. However, letting one of the parts of the boundary be the empty region of the boundary, which corresponds to the trivial word, we see that the word formed by the entire boundary must be equivalent to the trivial word.

\end{proof}

\section{Conjecture 1 for $A_n$} \label{secAn}

We start by fixing the following Coxeter-Dynkin diagram for $A_n$ for the rest of our paper.

\begin{figure}[H]
\centering
\includegraphics[scale = 0.5]{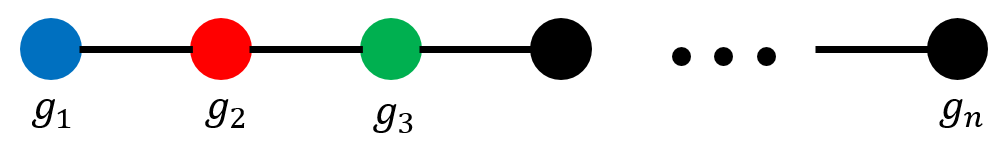} 
\caption{A colored Coxeter-Dynkin diagram for $A_n$. We note that in the above figure and throughout the rest of this paper, black will denote an arbitrary color that has not yet been assigned to a generator.}
\label{coloredAn}
\end{figure}

\begin{Lemma}\label{removefirstgen}
Let the Coxeter group $A_n$ have generators $g_1, \ldots, g_n$. Every word in $A_n$ can be written with at most 1 occurrence of $g_1$.
\end{Lemma}

\begin{proof}
This can be proven easily by using induction on the length of the word.
\end{proof}

\begin{Definition}
Given a pair of colors in our diagram, the two colors are said to \emph{commute} if their corresponding generators commute.
\end{Definition}

\begin{Definition}
Given $c$ an arbitrary color, we let the \emph{$c$ subgraph} of a diagram denote the graph on the vertices that have adjacent edges of color $c$ and all the $c$ edges in which we ignore vertices of degree 2 by gluing the edges.
\end{Definition}

\begin{Remark}
We see that all vertices of degree 4 with the color $c$ in our original diagram are vertices with colors $c/d$ where $d$ is a color that commutes with $c$. These vertices are vertices of degree 2 in the $c$ subgraph and are thus ignored. As a result, when considering the $c$ subgraph, we ignore vertices formed by colors that commute with $c$.

\end{Remark}

\begin{figure}[H]
\centering
\includegraphics[scale = 0.65]{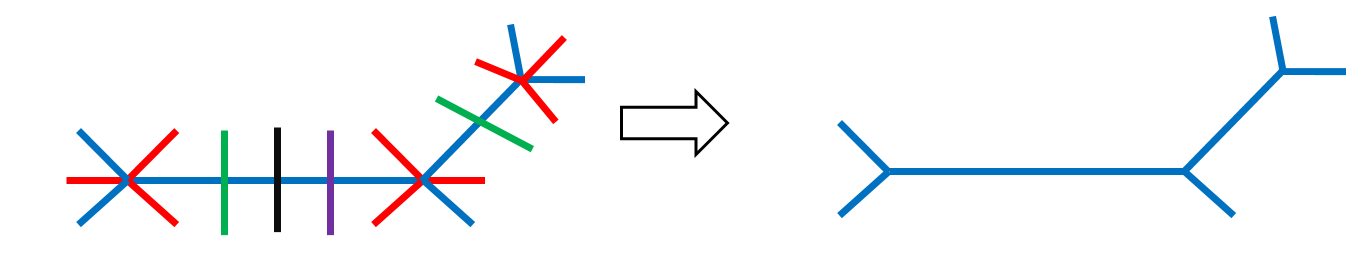}
\caption{The picture on the right is the blue subgraph of the subdiagram to the left.}
\label{comp}
\end{figure}

\begin{Example}
Figure~\ref{comp} provides an example of a subdiagram.

\end{Example}

\begin{Definition}

Letting $c$ be an arbitrary color, we let \emph{$c$-face} stand for a face in the $c$ subgraph.
\end{Definition}

\begin{Definition}
Two vertices of the same type with colors $c_1$ and $c_2$ are called \emph{almost $c_1$-adjacent} if they are connected by an edge $E$ in the $c_1$ subgraph.
\end{Definition}

\begin{Lemma}\label{almostadj}

Given two almost blue-adjacent vertices of type blue/red in a diagram for $A_n$, either they can be deleted or they can be transformed into the diagram in Figure~\ref{2blueverts}.

\begin{figure}[H]
\centering
\includegraphics[scale = 0.2]{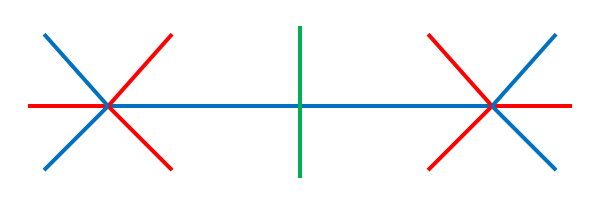}
\caption{2 almost-blue adjacent vertices of type blue/red.}
\label{2blueverts}
\end{figure}
\end{Lemma}

\begin{proof}
Let the blue edge connecting the 2 blue/red vertices in the blue subgraph be denoted as $E$. We see that all the vertices on the edge $E$ (in the context of the entire diagram and not just the blue subgraph) correspond to a word in $A_n$ comprised of generators $g_3, \ldots, g_n$. The group generated by $g_3, \ldots , g_n$ with relations between these generators as given in the  Coxeter-Dynkin diagram in Figure ~\ref{coloredAn} is isomorphic to the Coxeter group $A_{n-3}$. Using Lemma ~\ref{removefirstgen}, we can rewrite any word in $A_{n-3}$ with at most 1 instance of $g_3$. Thus we can replace the vertices on the edge $E$ with equivalent vertices, where there is at most 1 vertex of type blue/green. However, given any vertex of type blue/$c$ where $c$ is a color that is not green, we see that $c$ commutes with both red and blue, and thus can be moved out of the edge $E$ using the ZAM relation in Figure ~\ref{flipzam}. 

Thus, moving out all colors $c$ that commute with red and blue, we are left with at most 1 vertex on edge $E$ that is of type green/blue. If there are no vertices on edge $E$, then we have 2 adjacent vertices of the same type connected by an uninterrupted edge, and thus we can use Lemma ~\ref{adjvert} to remove them. In the other case, we are left with precisely the diagram in the statement of this lemma (Figure ~\ref{2blueverts}).
\end{proof}

\begin{Lemma}\label{trick}
The 2 diagrams in Figure \ref{2equivdiags} are equivalent.

\begin{figure}[H]\label{figtrick}
\centering
\includegraphics[scale = 0.3]{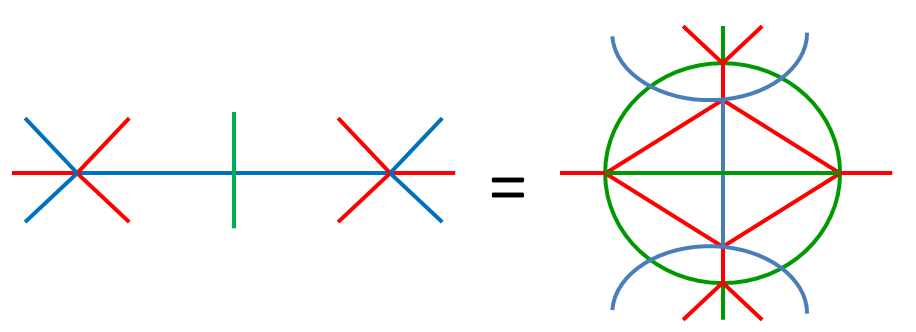}
\caption{2 equivalent diagrams.}
\label{2equivdiags}
\end{figure}

\end{Lemma}

\begin{proof}
To show that the two diagrams in Figure ~\ref{2equivdiags} are equivalent, we connect the boundaries of the 2 diagrams, and note that this is possible as both diagrams have the same boundary. After connecting the boundaries of the 2 diagrams, we get a closed diagram. Using the $A_3$ ZAM relation, we can show that the resulting diagram is equivalent to the trivial one. Using Lemma ~\ref{boundarycor}, we see that the 2 diagrams are equivalent.

\end{proof}
Having stated and proven the necessary lemmas, we now provide our proof for $A_n$ and how we solved certain challenges that arose. 
The fundamental idea of our proof is using induction on the number of blue/red vertices (vertices corresponding to generators $g_1$ and $g_2$). One of the major challenges we faced while trying to prove $A_n$ was that there are many different types of vertices. Initially, we tried to use ther Euler characteristic to find a small face and then show there must be vertices around this small face that we can delete. However, due to the many different types of vertices, this was impossible. We instead came up with Lemma ~\ref{trick}, which allows us to reduce the size of a face in the blue-graph. This strategy is \emph{very} powerful, as it allows us to consider only 1 case for a blue-face: a blue-face of size 2, from where we can find blue/red vertices to delete. Thus Lemma ~\ref{trick} is a very important tool that we have developed. This lemma can be applied to any Coxeter group that has a subgroup isomorphic to $A_n$ for some $n \geq 3$.  


\begin{Theorem}\label{An}
Conjecture ~\ref{main} holds for the family of Coxeter groups $A_n$.  
\end{Theorem}

\begin{remark}
The symmetric group which is one of the most common groups found throughout mathematics is a Coxeter group of type $A_n$. 
\end{remark}

\begin{proof} 

To prove that all diagrams for $A_n$ are homotopic to the empty diagram, we first induct on $n$. We see that our base case is $A_2$, which is isomorphic to $I_3$, which we have already proven from Corollary ~\ref{cordihedral}. Now we wish to show given an arbitrary diagram for $A_n$ that there is at least 1 blue/red vertex that we can delete. In the blue subgraph, consider a blue-face of arbitrary size. Let one vertex on the blue-face be $X$ and an adjacent vertex be $Y$. Also, let the blue edge between $X$ and $Y$ be $E$. In the conext of the entire diagram (and not just the blue subgraph) $X$ and $Y$ are almost-blue adjacent. By Lemma ~\ref{almostadj}, we can either delete $X$ and $Y$ in which case we are done, or we can transform $X$, $Y$ and $E$ to the diagram in Figure \ref{XYE}. 


\begin{figure}[H]
\centering
\includegraphics[scale = 0.3]{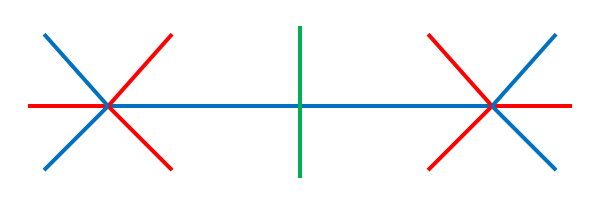} 
\caption{A diagram of $X$, $Y$ and $E$.}
\label{XYE}
\end{figure}

In the case where we have a blue/green vertex on $E$, we use Lemma ~\ref{trick} to reduce the size of the blue-face that we have been considering. It is not obvious how Lemma ~\ref{trick} actually reduces the size of the blue-face, thus we draw an example with a blue-face of size 4 in Figure ~\ref{trickfig}.

\begin{figure}[H]
\centering
\includegraphics[scale = 0.3]{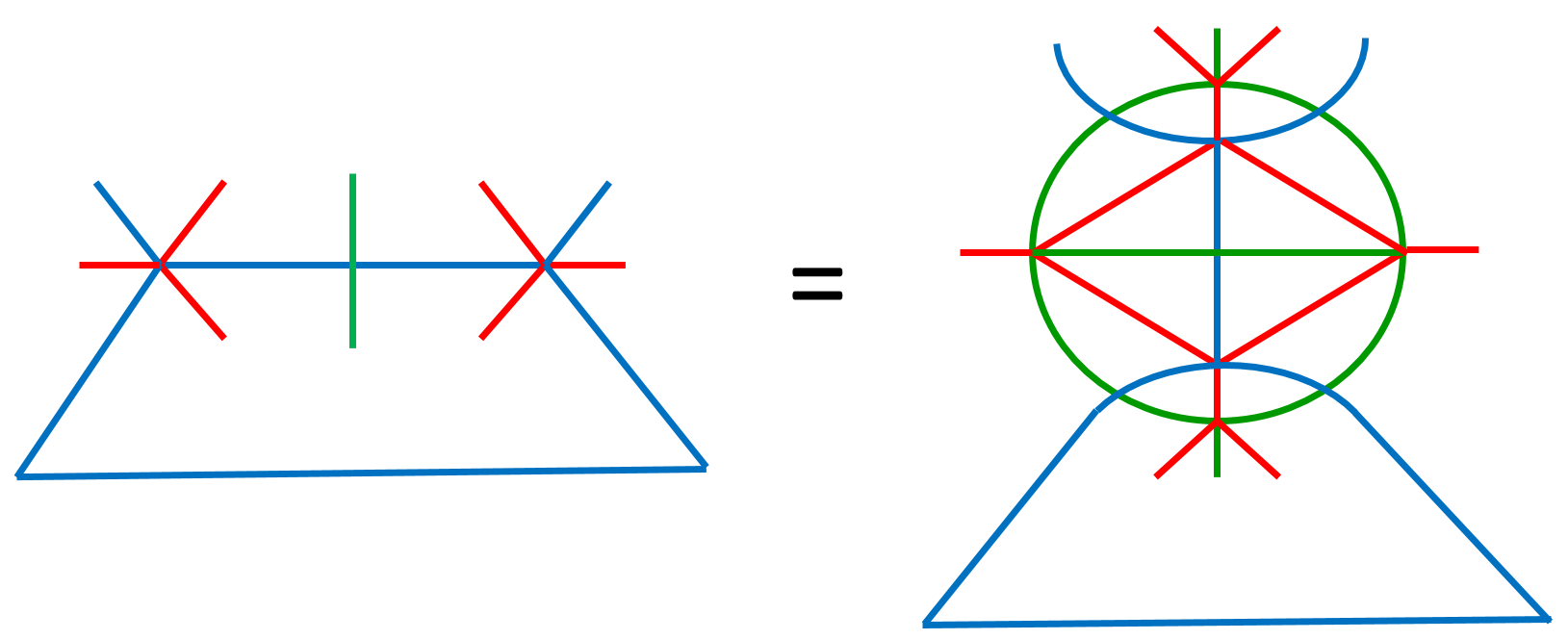} 
\caption{In the drawing above, we see that after applying Lemma ~\ref{trick}, we reduce the size of the blue-face from 4 to 3.}
\label{trickfig}
\end{figure}

We can continue this process of reducing the size of the blue-face until we are left with a blue-face of size 2. Using our method of removing colors that commute with both blue and red, there are only three possible blue-faces of size 2 where both red/blue cannot be removed, drawn below. We note that we have assigned to generator $g_4$ the color yellow, and see that for the group $A_3$, we do not have a generator $g_4$, and thus only one of the following blue-faces (the one with no yellow edges) is possible for the case of $A_3$.

\begin{figure}[H]\label{blue-faces_2}
\centering
\includegraphics[scale = 0.5]{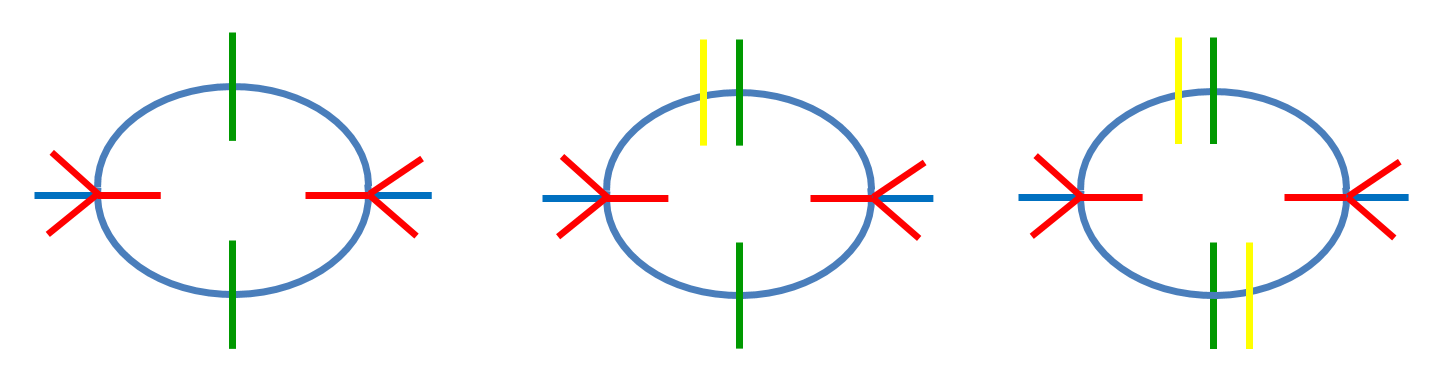}
\caption{The three possible blue-faces of size 2.}
\end{figure}

Recall that red, green and yellow represent $g_2, g_3$ and $g_4$ respectively. We know from Lemma ~\ref{trivialboundary} that the boundary of any subdiagram must be equivalent to the trivial word. However, it is simple to show that the boundaries of all of the blue-faces in Figure ~\ref{blue-faces_2} are not equivalent to the identity. Thus none of the 3 blue-faces in Figure ~\ref{blue-faces_2} are possible, and thus the only blue-face of size 2 that is possible is drawn below:

\begin{figure}[H]
\centering
\includegraphics[scale = 0.3]{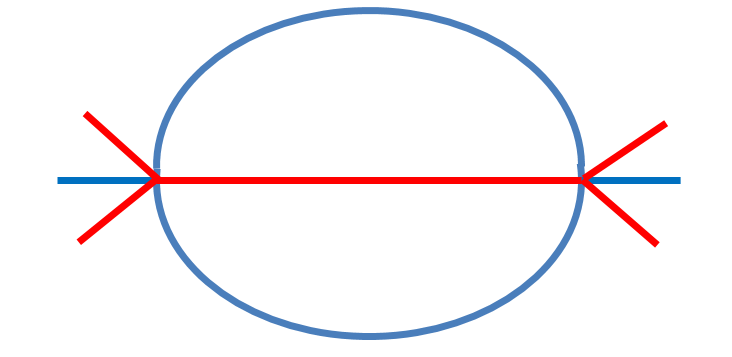}
\caption{The only possible blue-face of size 2.}
\label{threebluefaces}
\end{figure}

We note that in the above picture, we can assume the blue-face has nothing inside it except 2 connected red edges by the following reasoning: All of the colors inside the blue-face correspond to generators $g_2, \ldots, g_n$, which generate the group $A_{n-1}$. However, by our inductive hypothesis stated at the beginning of this proof, we can assume $A_{n-1}$ satisfies Conjecture ~\ref{main}, and thus by Corollary ~\ref{boundary}, all subdiagrams for $A_{n-1}$ with the same boundary are equivalent. Thus all subdiagrams of $A_{n-1}$ that have 2 red edges are equivalent to the one where the 2 red edges are connected, as in the picture above.

Since the blue-face in Figure ~\ref{threebluefaces} has adjacent vertices of the same type, by Lemma ~\ref{adjvert}, we can delete these blue/red vertices. Using induction on the number of blue/red vertices, in any diagram for $A_n$ we can delete all blue/red vertices until there are none left. After this, all the remaining blue vertices in a diagram are blue/$c$, where $c$ is a color that commutes with blue. Thus we can trivially remove all remaining blue vertices. Therefore, we have no more blue edges in our diagram, and our diagram is now a diagram for the group $A_{n-1}$. However now we can use our inductive hypothesis and reduce our diagram to the empty graph. Thus, all diagrams for the group $A_n$ are homotopic to the empty graph.

\end{proof}

\section{Conjecture 1 for $\I_n$} \label{In}

Finally, we prove Conjecture ~\ref{main} for the dihedral Artin braid groups. While this group has only two generators, it is a very difficult case due to the orientations on our planar graphs. The main challenge we faced is that two adjacent vertices of the same type are not necessarily deletable because orientations of the vertices may not line up. Our key insight in this case was to look at angles in our graph and rephrase our problem in terms of properties of these angles. These same properties of angles that we examine have been useful in our work thus far in other oriented cases. 

Recall for $\I_n$ we assign $g_1$ the color blue and $g_2$ the color green, as in Definition ~\ref{braidIn}.

\begin{Definition}
An angle of a vertex in a diagram is called \emph{varied} if the edges of the angle have different orientations. An angle is called \emph{uniform} if the edges of the angle have the same orientation.
\end{Definition}

\begin{figure}[H]
\centering
\includegraphics[scale = 0.42]{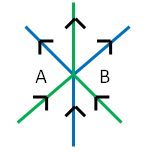} 
\caption{In this vertex for $\I_3$, angles $A$ and $B$ are varied angles. The rest of the angles are uniform angles.}\label{variedref}
\end{figure}

\begin{Example}
Figure~\ref{variedref} provides an example of varied and uniform angles. 
\end{Example}

\begin{Theorem}\label{braidInproof}
Conjecture ~\ref{main} holds for the braid groups $\I_n$.
\end{Theorem}
\begin{proof} 
We note that if a face in our graph has 2 adjacent varied angles, then we can use the bridge relation to obtain a face of size 2, which we can then remove using the $\I_n$ cancellation of pairs relation. Also, we see that a face of size 2 has either 2 varied angles (in which case we can delete it) or 2 uniform angles. Having made these observations, we prove that in any given diagram, there is at least 1 vertex that we can delete.

We first let $F$ stand for the number of faces in our graph, $E$ stand for the number of edges and $V$ stand for the number of vertices. Assume for the sake of contradiction that in a given diagram there are no vertices we can delete. We see that every face must have no 2 adjacent angles be varied angles. Thus all faces must have at least 2 uniform angles (including faces of size 2, since faces of size 2 have either 2 varied or 2 uniform angles, and if they have 2 varied angles the vertices of the face can be deleted using the directed cancellation of pairs relation). Thus the number of uniform angles is $\geq 2F$. However using the Euler formula for planar graphs, we know that $V+F=E+2$, where in this case $E=\frac{2nV}{2}=nV$, since each vertex of our graph has degree $2n$. Thus $F= (n-1)V+2$. Thus the number of uniform angles is $\geq 2F = 2((n-1)V+2) = (2n-2)V+4$. But every vertex has precisely $2n-2$ uniform angles, thus the number of uniform angles is exactly $(2n-2)V$. Thus the number of uniform angles is simultaneously $\geq (2n-2)V+4$ but equal to $(2n-2)V$. Clearly, this is a contradiction, and thus these exists a face with 2 adjacent varied angles, and we have found 2 vertices we can delete (using the $\I_n$ cancellation of pairs relation). Using induction, we can delete all vertices from our graph. Therefore, every diagram for $\I_n$ is homotopic to the empty diagram.

\end{proof}

\section{Conclusion and Future Directions}\label{conc}

The work presented in this paper uses the diagrammatics of Coxeter groups and braid groups to diagrammatically prove the $K(\pi,1)$ conjecturette for several families of Coxeter groups: $I_n$ and $A_n$. Additionally, we diagrammatically prove this conjecture for the braid group $\I_m$. This work addresses a question posed in \cite{elias} regarding the existence of a diagrammatic proof of the $K(\pi,1)$ conjecturette, which the authors of \cite{elias} were unable to find an elementary proof for. Perhaps the most important aspect of our work is the development of lemmas and strategies that can be used to tackle further cases of the $K(\pi,1)$ conjecturette for other Coxeter groups and braid groups---especially for cases that remain unsolved with traditional approaches. In particular, Lemma ~\ref{trick} can be applied to any Coxeter group that has a subgroup isomorphic to $A_n$ for $n \geq 3$. Additionally, our approach towards the oriented case $\I_m$ has yielded partial results for other braid groups. With this in mind, one might consider generalizing our work to the type $D$ Coxeter groups or extending our methods to other families of braid groups.

\section{Acknowledgements}
The authors would like to thank their mentor Alisa Knizel for her valuable insight and guidance. They also wish to thank the MIT PRIMES program for facilitating this research opportunity, and in particular wish to acknowledge Dr. Khovanova, Dr. Gerovitch and Professor Etingof. Finally, the authors thank Professor Ben Elias for suggesting this project and providing helpful guidance.

\newpage

\end{document}